\documentclass{article}

\usepackage{latexsym}
\usepackage{amsmath}
\usepackage{verbatim}
\usepackage{amssymb}
\usepackage{amsthm}
\usepackage{graphicx}

\usepackage[top=3cm, bottom=3cm, left=2cm, right=2cm]{geometry}

\numberwithin{equation}{section}

\newcommand{\abs}[1]{\left|#1\right|} 
\newcommand{\R}{{\mathbb {R}}}
\newcommand{\norm}[2]{\left\| #1 \right\|_{#2}}

\newcommand{\bC}{\mathbf{C}}

\newtheorem{theorem}{Theorem}
\newtheorem{proposition}{Proposition}[section]

\newtheorem{lemma}[proposition]{Lemma}
\newtheorem{remark}[proposition]{\bf Remark}
 
\DeclareMathOperator{\sgn}{sgn}

\title{ The one-dimensional Keller-Segel model with fractional diffusion of cells\footnote{During the final stages of the preparation of this paper, the authors were 
made aware of similar works by Biler \& Laurencot, and also Rodrigo.
}}

\author{Nikolaos Bournaveas\protect\footnote{University of Edinburgh, School of Mathematics,
JCMB, King's Buildings, Edinburgh EH9 3JZ, UK.
 E-mail address:
\texttt{n.bournaveas@ed.ac.uk}}, Vincent Calvez\protect\footnote{Unit\'e de Math\'ematiques Pures et Appliqu\'ees, \'Ecole Normale Sup\'erieure de Lyon, 
46 all\'ee d'Italie
69007 Lyon,
France. E-mail address: \texttt{vincent.calvez@ens-lyon.fr}}\;\footnote{Corresponding author}}

\date{\today, {\em in progress}}

\begin{document}

\maketitle

\begin{abstract}
We investigate the one-dimensional Keller-Segel model where the diffusion is replaced by a non-local operator, namely the fractional diffusion with exponent $0<\alpha\leq 2$. We prove some features related to the classical two-dimensional Keller-Segel system: blow-up may or may not occur depending on the initial data. More precisely a singularity appears in finite time when $\alpha<1$ and the initial configuration of cells is sufficiently concentrated. On the opposite, global existence holds true for $\alpha\leq1$ if the initial density is small enough in the sense of the $L^{1/\alpha}$ norm. 
\end{abstract}
 
\noindent{{\bf Keywords.} Self-organization, chemotaxis, fractional diffusion, global existence, blow-up.}

\section{Introduction}
Chemotaxis is the directed motion of cells in response to various chemical clues. It plays a key role in developmental biology, and more generally in self-organization of cell populations.
Several categories of mathematical models have been proposed to describe this organization process. Depending upon the level of description required, micro-, meso- or macroscopic models can be used \cite{Perthame04,Hierarchy,PerthameBook}. Mesoscopic models consist in kinetic (scattering) equations well-suited for describing the motion of bacteria such as {\em Escherichia coli} which undergo a run and tumble process \cite{ErbanOthmer04,BC}. Macroscopic models consist in parabolic (drift-diffusion) equations and are well-suited for describing motion of large cells such as the slime mold amoebae {\em Dictyostelium discoideum} \cite{Horstmann,ErbanOthmer07,HillenPainter09}. We focus on the macroscopic setting in this paper.

The so-called Keller-Segel model exhibits a very rich behaviour, emphasized by the critical mass phenomenon arising in two dimensions of space. The Keller-Segel system writes in a simple formulation \cite{JL}:
\begin{subequations}\label{PKS}
\begin{align}
 \partial_t\rho(t,x)& =  \Delta\rho(t,x) - \nabla \cdot \left( \rho(t,x) \nabla c(t,x) \right)\, ,\quad t>0\, , \quad x\in \R^d \label{PKSa}\\
- \Delta c(t,x) & =  \rho(t,x) \, .\label{PKSb}
\end{align}
\end{subequations}
Here $\rho(t,x)$ denotes the cell density and $c(t,x)$ denotes the concentration of the chemical attractant. The first contribution in the right hand side of \eqref{PKSa} expresses the tendency of cells
to diffuse under their own Brownian motion whereas the second term expresses their tendency to aggregate due to the
presence of the chemical. In two dimensions of space the two tendencies are evenly balanced
and the global behavior of the solution depends on the total mass of the cells. More precisely, for
$M> 8 \pi$ blow-up occurs in finite time (aggregation overwhelms diffusion) and for $M<8\pi$ solutions are global in time (diffusion wins the competition) \cite{BDP}. 

However in one  dimension of space diffusion is always stronger than 
aggregation and blow-up never occurs for systems such \eqref{PKS} \cite{HP, OY}.
 

In this paper we study the system \eqref{PKS} in one space dimension 
with the cell diffusion being ruled by fractional diffusion. The usual Laplacian in \eqref{PKSa} is therefore replaced by the fractional Laplacian. The non-local parabolic equation writes as following:
\begin{subequations}\label{frKS}
\begin{align}
\partial_t \rho(t,x) & = - \Lambda^{\alpha} \rho(t,x) - \partial_x \left( \rho \partial_x c  \right)\, , \quad t>0\, , \quad x\in \R \label{eqrho}\\
- \partial^2_{xx} c(t,x) & =  \rho(t,x) \, , \label{eqc} 
\end{align}
\end{subequations}
equipped with suitable initial condition $\rho(0,\cdot) = \rho_0$ and decay conditions at infinity.
For an exponent $\alpha\in (0,2]$, the positive operator $\Lambda^\alpha = (-\Delta)^{\alpha/2}$ is defined in Fourier variables by
$\widehat{\Lambda^\alpha f}(\xi)=|\xi|^\alpha \widehat{f}(\xi)$. An alternative representation is given by: 
\[ \Lambda^\alpha f(x) = c_\alpha \int_{y\in \R} \dfrac{f(x)-f(y)}{|x-y|^{1+\alpha}}\, dy = {c_\alpha} \int_{h\in \R} \dfrac{2f(x) - f(x+h)-f(x-h)}{|h|^{1+\alpha}} \, dh\, , \]
where $c_\alpha$ is some normalizing factor.

Non-local operators, and in particular the fractional Laplacian, have received
a lot of attention  recently \cite{CSi, CSou, CV}. In biology the motivation comes from the fact 
that in many cases organisms adopt L\'evy-flight search strategies and therefore 
dispersal is better modelled by non-local operators  \cite{B, E1, E2,  K, L}. Focusing on the one-dimensional case may seem 
unnatural
 from the biological viewpoint. However we have in mind seeking a critical mass phenomenon as it has been derived for the 
 two-dimensional classical Keller-Segel model \eqref{PKS}. It appears that $\alpha = d$ is the critical fractional exponent to 
 state such a result. Therefore it makes only sense when $d=2$ or $d=1$. We partially answer this issue below in the latter 
 situation.

The system \eqref{frKS} was first studied in \cite{E1} where it was shown that global
existence holds true for $1<\alpha \leq 2$ assuming that $\rho_0 \in L^1 \cap L^2$ and $\rho_{0}'\in L^2$.

We aim at providing here global existence versus blow-up results in the same spirit as for the dichotomy arising in the two-dimensional classical Keller-Segel system \eqref{PKS}. More precisely, we are able to prove that solutions are global in time in the 'fair-competition' case $\alpha = 1$, if  the total mass $M$ is assumed to be small enough. In the case $\alpha<1$ we show that solutions may exist globally or may blow-up depending on the initial data. We exhibit explicit criteria to distinguish between these two alternatives. For the case $\alpha>1$ we improve previous work \cite{E1} by weakening regularity hypotheses on the initial data. 




\begin{theorem}[Global existence]\label{global}
Consider the system \eqref{frKS} for $0< \alpha\leq 1$ 
with initial data $\rho_0\in L^{p_0}(\R)$ for some $p_0>1/\alpha$. There 
exists a constant $K_1(\alpha)$ such that the condition, \[\|\rho_0\|_{L^{1/\alpha}}<K_1(\alpha)\, ,\]  guarantees existence of global weak solutions.\\ In addition, regularizing effects act for \eqref{frKS}, and the density belongs to any $L^p$ space for any positive time $t>0$. \\
In the case $1< \alpha\leq 2$, assume $\rho_0\in L^{p_0}(\R)$ for some $p_0>1$. Then solutions are global in time and belong to any $L^p$ space for all positive time $t>0$. 
\end{theorem}

To complete the picture it is natural to look for blow-up results in the super-critical case. We shall prove in the sequel that the 
aggregation contribution can overcome the diffusion effect in the case $\alpha<1$ under suitable restrictions on the initial data. 
However describing the behaviour for initial data having large mass in the case $\alpha = 1$ remains open (the constant $K_2(\alpha)$ in Theorem \ref{BU} below diverges when $\alpha\to 1$).

\begin{theorem}[Blow-up]
\label{BU}
Consider the system \eqref{frKS} for $0< \alpha< 1$ in one space dimension 
with initial data $\rho_0\in L^1((1+|x|)dx)$. There exists a constant $K_2(\alpha)$ such that the condition,
\[ \left( \int_\R |x|\rho_0(x)\, dx\right)^{1-\alpha} \leq K_2(\alpha) M^{2-\alpha}\, , \]
excludes global existence of regular solution: a singularity must appear in finite time. 
\end{theorem}

The paper is organized as follows: in Section \ref{sec:GE} we prove global existence, beginning with a simple but not complete 
argument based on $L^2$ estimates. The proof is then achieved thanks to $L^p$ estimates inspired by \cite{CV}. 
In Section \ref{sec:BU} we prove blow-up of solutions. The paper is supplemented by numerical illustrations of the two 
above-mentioned phenomena.



\thanks{{\em Acknowledgements.} 
NB would like to thank the Laboratoire Jacques-Louis Lions of the Universit\'e Pierre et Marie Curie (Paris, France)
and the D\'epartement de Math\'ematiques et Applications of the \'Ecole Normale Sup\'erieure (Paris, France) for their 
hospitality and financial support during his sabbatical leave in the spring semester of 2008, during which 
part of the research for this paper was done. VC thanks the Centre de Recerca Matem\'atica (Bellaterra, Spain), for providing an excellent atmosphere of research during the research program 'Mathematical Biology'.}

\section{Global existence for small initial data: proof of Theorem \ref{global}}

\label{sec:GE}

We start by stating some estimate which will be widely used throughout this section.
\begin{proposition}[Interpolation inequality] \label{prop:GNS}
For any exponents $0< \alpha\leq 1$ and $1\leq p < +\infty $, the following Gagliardo-Nirenberg type inequality holds true:
\begin{equation}
\int_\R \rho^{p+1}(x)\, dx  \leq \bC(p,\alpha) \norm{\rho^{p/2}}{\dot H^{\alpha/2}}^2 \norm{\rho}{1/\alpha}\, .
\end{equation}
\end{proposition}

\begin{proof}
We distinguish between the cases $\alpha<1$ and $\alpha = 1$. In the former we use first the H\"older inequality to obtain:
\begin{align*}
\int_\R \rho^{p+1}(x)\, dx & \leq \left(\int_\R  \rho^{p /(1-\alpha)}(x)\, dx \right)^{1-\alpha} \left(\int_\R  \rho^{1/\alpha}(x)\, dx \right)^{\alpha}\\
& \leq \bC(p,\alpha) \norm{\rho^{p/2}}{\dot H^{\alpha/2}}^2 \norm{\rho}{1/\alpha}\, ,
\end{align*}
where we have used the Sobolev embedding: $\dot H^{\alpha/2}\hookrightarrow L^{2/(1-\alpha)}$.

In the case $\alpha =1$, we can use the following general result 
\cite{mo}: for any $\lambda, \mu, s, q, r, \theta \in \R$ satisfying the following  relations:
 \[
  1\leq s,q\leq r \leq \infty\ ,\ 0<\theta<1\ ,\ \lambda>\frac{d}{s} - \frac{d}{r}\ ,
\ \mu < \frac{d}{q} - \frac{d}{r} ,
  \]
  \[
 \theta \left(\lambda - \frac{d}{s} + \frac{d}{r} \right) + 
(1-\theta) \left(\mu - \frac{d}{q} + \frac{d}{r}  \right)=0 .
  \]
we have,
\begin{equation}\label{gn6}
\norm{f}{L^r} \leq \bC  \norm{f}{\dot{W}^{\lambda,s}}^{\theta} \norm{f}{\dot{W}^{\mu,q}}^{1-\theta} .
\end{equation}
Applying that to the particular choice: $f = \rho^{p/2}$, $\lambda = 1/2$, $\mu=0$, $s = 2$, $r = 2(p+1)/p$, $\theta = 2/r$, $q=2/ p$ yields the result.

Observe that proceeding as above, the exponent $p$ cannot be chosen arbitrarily (the constraint $q\geq1$ forces $p\leq 2$). 
However, there is a way to extend it to any $p\geq 1$ by slightly modifying the argument: using $f = \rho^{p/2}$, $\lambda = 1/2$, $\mu=0$, $s = 2$, $r = 2(p+1)/p$, $\theta = 1/(p+1)$, $q=2$ we get:
\begin{align*}
\left( \int_\R \rho^{p+1}(x)\, dx \right)^{p/(2(p+1))} & \leq C \norm{\rho^{p/2}}{\dot H^{1/2}}^{1/(p+1)} \left( \int_\R \rho^{p}(x)\, dx \right)^{p/2(p+1)} \\
& \leq C \norm{\rho^{p/2}}{\dot H^{1/2}}^{1/(p+1)} \left(\norm{\rho}{1} \left(\int_\R \rho^{p+1}(x)\, dx\right)^{p-1}\right)^{1/2(p+1)}\, .
\end{align*}
Raising this inequality to the power $2(p+1)$ leads to the result.
\end{proof}

\subsection{{\em A priori} $L^2$ estimates}\label{HE}

We complete here some existing results first derived by Escudero  \cite{E1}. We use Gagliardo-Nirenberg type
inequalities instead of the Sobolev inequality used in \cite{E1}. This allows us to study a wider range of $\alpha$'s.
We are concerned in this section with the global existence of the Keller-Segel system with fractional diffusion of cells when $1/2\leq\alpha\leq 1$, 
using simple harmonic analysis estimates. This will be extended below in Section \ref{sec:existence Lp} to any $0<\alpha\leq 1$. The 
purpose of this section is to derive simply {\em a priori} estimates which guarantee global existence of solutions and to set the 
stage for our approach in Section \ref{sec:existence Lp}. The constraints on the exponent $\alpha$ here are an artefact of the 
method: 
in short the interpolation of $L^{1/\alpha}$ between $L^1$ and $L^2$ yields $1/2 \leq \alpha$.

As it is now standard in such systems, we aim at deriving suitable $L^p$ norm of the cell density. Due to the simple formulation of 
the fractional diffusion in the Fourier space variable, we opt for $p = 2$. We will relax this constraint in the next section.
We have the following estimation:
\begin{align*}\label{E}
 \frac{d}{dt}\dfrac12\norm{\rho(t)}{L^2}^{2}  & =  \int_{\R} \left(- \Lambda^{\alpha} \rho(t,x) - \partial_x \left( \rho(t,x) \partial_x c(t,x)  \right)\right) \rho(t,x)\, dx \\
  & = -  \int_{\R} \left(\Lambda ^{\alpha/2} \rho(t,x)\right)^2\, dx  + \frac12 \int_\R   \rho ^3(t,x)\, dx\, .
\end{align*}
We then apply the Gagliardo-Nirenberg inequality (Proposition \ref{prop:GNS}) for $p=2$: 
\begin{equation}\label{eq:GNS applied}
 \int_\R \rho(t,x)^3 \, dx \leq \bC(2,\alpha)   
 \norm{\Lambda^{\alpha/2} \rho(t)}{L^2}^{2} \left(\int_{\R} \rho^{1/\alpha}(t,x)\, dx\right)^\alpha\, .
\end{equation}
In the case $\alpha=1$ we obtain the decay of the $L^2$ norm providing that the mass is small enough:
\begin{equation}\label{eq:harmonic JL}
 \frac{d}{dt}\dfrac12\norm{\rho(t)}{L^2}^{2} \leq \left(-\dfrac1{\bC(2,1)M} +\dfrac12 \right) \int_\R \rho^3(t,x)\,  dx\, .
\end{equation}
It follows that $\norm{\rho(t)}{L^2}\leq \norm{\rho_0}{L^2}$, as soon as $\rho_0\in L^2$. It is also possible to conclude without assuming $\rho_0\in L^2$, by means of regularizing effects. In fact using interpolation between $L^1$ and $L^3$ it comes out that \eqref{eq:harmonic JL} also implies (when the mass is small enough):
\[ 
 \frac{d}{dt}\dfrac12\norm{\rho(t)}{L^2}^{2} \leq \left(-\dfrac1{\bC(2,1)M} +\dfrac12 \right) M^{-1} \norm{\rho(t)}{L^2}^4\, .
  \]
Therefore $\|\rho(t)\|_2$ becomes finite in zero time. We shall come back to that later.

In the case $\alpha <1$ the Gagliardo-Nirenberg inequality \eqref{eq:GNS applied} implies that:
\[ \dfrac d{dt } \dfrac12\norm{\rho(t)}{L^2}^{2} \leq \left( - \dfrac{1}{ \bC(2,\alpha) \|\rho(t)\|_{L^{1/\alpha}}} + \dfrac12 \right) \int_\R \rho^3(t,x)\, dx\, . \] 
As opposed to the case $\alpha = 1$, the quantity $\|\rho(t)\|_{L^{1/\alpha}}$ is not conserved in time. Therefore we have to develop an alternative strategy as in \cite{CPZ} for the Keller-Segel in dimension $d>2$, where the criterion for global existence involves the $L^{d/2}$-norm. Here we simply use the fact that $L^{1/\alpha}$ can be interpolated between $L^1$ and $L^2$ if $1/2\leq \alpha\leq 1 $. As a consequence we have:
\[ \dfrac d{dt }\dfrac12 \norm{\rho(t)}{L^2}^{2} \leq \left( - \dfrac{1}{\bC(2,\alpha) M^{2\alpha - 1} \|\rho(t)\|_{L^{2}}^{ 2-  2\alpha}} + \dfrac12 \right) \int_\R \rho^3(t,x)\, dx\, .  \]
Thus if the quantity $M^{2\alpha - 1} \|\rho_0\|_{L^{2}}^{ 2-  2\alpha}$ is small enough, then $\|\rho(t)\|_{L^2}$ automatically decays for every time. We will see later that this criterion can be ameliorated, as the $L^{1/\alpha}$ (before interpolation) appears to be the critical space for this problem (analogous to $L^{d/2}$ in the classical Keller-Segel problem). To derive this improved criterion we shall understand how the $L^p$ norms of the cell density evolve, using more refined tools for integration by parts.

\begin{remark}\label{rem}
In the case $1\leq \alpha \leq 2$, if we assume that $\rho_0 \in L^1 \cap L^2$ and $\rho_{0}'\in L^2$ we can work similarly as in \cite{E1}
to obtain an a-priori estimate on  $\norm{\rho_{x}(t)}{L^2(\R)}$. Then the Sobolev inequality gives a bound on 
$\norm{\rho(t)}{L^\infty}$.  
\end{remark}

\subsection{{\em A priori} $L^p$ estimates}

\label{sec:existence Lp}

Following \cite{CSi, CV}, the one-dimensional fractional Laplacian can be interpreted as a `Dirichlet to Neumann problem' on the 
two-dimensional half-space (with an appropriate modification when $\alpha\neq 1$). Namely it is related to the following 
minimization problem. Given a function $\rho(x)$ defined for $x\in \R$ (and belonging to appropriate spaces, see \cite{CSi} for 
details) find a function $\rho_*(x,y)$ defined on $\R\times (0,\infty)$ coinciding with $\rho(x)$ on the boundary: 
$\rho_*(x,0) = \rho(x)$,
which minimizes the weighted functional, 
\[ J(u) = \dfrac12\int_0^\infty\int_\R |\nabla u(x,y)|^2 y^{1-\alpha} \, dx dy\, . \]
When $\alpha = 1$ this is nothing but the harmonic extension of $\rho$ to the upper half-space. The fractional Laplacian is then deduced from the normal derivative of $\rho_*(x,y)$ on the boundary $\{y=0\}$ as described below. We will strongly use this minimization property. 

\begin{proposition}[Integration by parts: fractional diffusion \cite{CV}]
\label{prop:IPP}
Assume $\rho(x)$ is regular, then the following estimate holds true:
\[
\int_\R \rho^{p-1}(x) \Lambda^\alpha\rho(x)\, dx     \geq  \dfrac{4(p-1)}{p^2} \left\| \rho^{p/2}\right\|^2_{\dot H^{\alpha/2}}\, . \]
\end{proposition}

\begin{proof}
For the sake of completeness, we recall the main lines of the proof of Proposition \ref{prop:IPP}. We begin with the case $\alpha = 1$ which is somewhat simpler.
\paragraph{The half-Laplacian.}
In short, the one-dimensional half-Laplacian $\Lambda \rho$ is the normal derivative of the harmonic extension on the upper-half plane of $\rho$:
\begin{align*}
&\Lambda \rho(x) = - \partial_y \rho_* (x,0)  \, , \\
\mbox{where}\quad &\left\{\begin{array}{l} -\Delta \rho_*(x,y) = 0 \quad \mbox{on}\quad \R\times (0,\infty)\, ,\medskip \\
\rho_*(x,0) = \rho(x)\, .\end{array}\right.
\end{align*}

Using this characterization, we are able to integrate by parts and to estimate the following diffusion contribution (which appears in the proof of Theorem \ref{global} below): 
\begin{align}
 \int_\R \rho^{p-1}(x)\Lambda \rho(x)\, dx  & =   \int_\R \rho_*^{p-1}(x,0) \nabla \rho_*(x,0)\cdot \nu \, dx  \nonumber\\
& =  \int_0^\infty\int_\R \nabla \rho_*^{p-1}(x,y)\cdot\nabla \rho_*(x,y)\, dxdy  \nonumber \\
& =  \dfrac{4(p-1)}{p^2} \int_0^\infty\int_\R |\nabla \rho_*^{p/2}(x,y)|^2\, dxdy\nonumber \\
& \geq \dfrac{4(p-1)}{p^2}  \int_0^\infty\int_\R |\nabla (\rho^{p/2})_*(x,y)|^2\, dxdy \nonumber\\  
&  \geq \dfrac{4(p-1)}{p^2} \left\|\rho^{p/2}\right\|^2_{\dot H^{1/2}}\, . \label{eq:1/2 IPP}
\end{align}

\paragraph{The $\alpha/2-$Laplacian.}
For any $0<\alpha<2$ the fractional Laplacian $\Lambda^\alpha\rho$ can be interpreted as follows \cite{CSi}: 
\begin{align*}
&\Lambda^\alpha \rho(x) = \lim_{y\to 0} \left[- y^{1-\alpha}\partial_y \rho_* (x,y)\right]  \, , \\
\mbox{where}\quad &\left\{\begin{array}{l} -\nabla\cdot\left(y^{1-\alpha} \nabla\rho_*\right) (x,y) = 0 \quad \mbox{on}\quad \R\times (0,\infty)\, , \medskip\\
\rho_*(x,0) = \rho(x)\, .\end{array}\right.
\end{align*}

In the same lines as \eqref{eq:1/2 IPP} we are able to estimate the following diffusion contribution:
\begin{align}
\int_\R \rho^{p-1}(x) \Lambda^\alpha\rho(x)\, dx    & =  \int_\R \rho_*^{p-1}(x,0) y^{1-\alpha} \nabla \rho_*(x,0)\cdot \nu \, dx  \nonumber \\
& =  \int_0^\infty\int_\R \nabla \rho_*^{p-1}(x,y)  \cdot y^{1-\alpha}\nabla \rho_*(x,y) \, dxdy \nonumber \\
& =  \dfrac{4(p-1)}{p^2}\int_0^\infty \int_\R |\nabla \rho_*^{p/2}(x,y)|^2  y^{1-\alpha} \, dxdy \nonumber \\
& \geq \dfrac{4(p-1)}{p^2} \int_0^\infty \int_\R |\nabla (\rho^{p/2})_*(x,y)|^2  y^{1-\alpha} \, dxdy \nonumber \\ 
&  \geq  \dfrac{4(p-1)}{p^2} \left\| \rho^{p/2}\right\|^2_{\dot H^{\alpha/2}}\, . \label{eq:alpha IPP}
\end{align}
\end{proof}

\begin{proof}[Proof of Theorem \ref{global}]
The case  $1<\alpha\leq 2$ has already been  treated in \cite{E1} (see Remark \ref{rem:alpha>1} at the end of the proof), so we focus on the situation where  $0<\alpha \leq 1$ in the sequel.

\paragraph{$L^{1/\alpha}$ is the critical space.} Following the lines of \cite{JL} and \cite{CPZ} we estimate the evolution of the $L^p$ norms of the cell density:
\begin{align*} 
\dfrac d{dt} \dfrac1p\| \rho(t) \|^p_{L^p}  &  = -  \int_{\R} \rho^{p-1}(t,x) \Lambda ^{\alpha} \rho(t,x)\, dx  +  
\frac{p-1}{p} \int_\R   \rho^{p+1}(t,x)\, dx \\
& \leq - \dfrac{4(p-1)}{p^2} \left\| \rho^{p/2}(t)\right\|^2_{\dot H^{\alpha/2}} + \frac{p-1}{p} \int_\R   \rho^{p+1}(t,x)\, dx\, .
\end{align*}

Using Proposition \ref{prop:GNS} we obtain:
\begin{equation}
 \int \rho^{p+1} (t,x) dx \leq \mathbf{C}(p,\alpha) \norm{\rho(t)}{L^{1/\alpha}} \norm{\rho^{p/2}(t)}{\Dot{H}^{\alpha/2}}^{2}\,  ,
\end{equation}
therefore
\begin{align}
\dfrac d{dt} \dfrac1p\| \rho(t) \|^p_{L^p} & \leq \left( - \dfrac{4(p-1)}{ p^2\bC(p,\alpha) \|\rho(t)\|_{L^{1/\alpha}}} + 
\frac{p-1}{p} \right) \int_\R \rho^{p+1}(t,x)\, dx\, .
\label{eq:Lp decays}
\end{align}
 
   Choosing in particular $p = 1/\alpha$ we 
obtain that the $L^{1/\alpha}$ norm is time-decreasing whenever $\| \rho_0 \|_{L^{1/\alpha}}$ is strictly smaller than 
$\frac{4\alpha}{\mathbf{C}(1/\alpha,\alpha)}$.

\paragraph{Regularizing effects.}  We shall prove within the next lines that the cell density $\rho(t,\cdot)$ belongs to any $L^p$ space for arbitrary positive time, provided that the initial $L^{p_0}$ norm is finite for some $p_0>1/\alpha$. The argument follows the main lines of \cite{JL,CPZ,CV}. 

First we shall relax the criterion on $\| \rho_0 \|_{L^{1/\alpha}}$ to 
\begin{equation}
\| \rho_0 \|_{L^{1/\alpha}} < \dfrac4{ p_0 \bC(p_0,\alpha)}\, .
\end{equation}
This ensures that the $L^{p_0}$-norm, which is initially finite by assumption, is decreasing in time. As a consequence, we get the following upper-bound for any truncation $k>0$:
\begin{equation} 
\|(\rho(t) - k)_+\|_{L^{1/\alpha}} \leq |\{x:\rho(t,x) > k\} |^{\alpha - 1/p_0}\,\|(\rho(t) - k)_+\|_{L^{p_0}}  \leq \left(\dfrac M{k}\right)^{\alpha - 1/p_0} \|\rho_0\|_{L^{p_0}}\, .
\label{eq:equi-integrability}
\end{equation}

Second, we extend the above strategy to the derivation of $\|(\rho(t) - k(p))_+\|_{L^p}$ for some $k(p)>0$ to be chosen later:
\[\begin{split} \dfrac d{dt} \dfrac1p\| (\rho(t) - k(p))_+ \|^p_{L^p}    \leq - \dfrac{4(p-1)}{p^2} \|(\rho(t) - k(p))_+\|_{\dot H^{\alpha/2}} + \dfrac{p-1}p \int_\R   (\rho(t,x) - k(p))^{p+1}\, dx \\ +  C(k,p) \int_\R   (\rho(t,x) - k(p))^{p}\, dx + C(k,p) \int_\R   (\rho(t,x) - k(p))^{p-1}\, dx   \, .
\end{split}\]
The last term can be interpolated between $L^1$ and $L^p$. The nonlinear contribution of homogeneity $p+1$ goes as previously, except that we shall ensure here that $ \|(\rho(t) - k(p))_+\|_{L^{1/\alpha}}$ is strictly smaller than $4/(p\bC(p,\alpha))$ independently of time. 

Introduce the notation: $Y_p(t) = \| (\rho(t) - k(p))_+ \|_p^p$. We have,
\[ \dfrac d{dt} Y_p(t) \leq \left(  - \dfrac{4(p-1)}{ p^2\bC(p,\alpha) \|(\rho(t) - k(p))_+\|_{L^{1/\alpha}}} + \dfrac{p-1}p 
\right) \int_\R (\rho(t,x) - k(p))_+^{p+1}\, dx + O\left(Y_p(t)\right) + O(1)\, . \]
Using the following interpolation inequality:
\[ Y_p(t)\leq M ^{1/p} \left(\int_\R (\rho(t,x) - k(p))_+^{p+1}\, dx\right)^{1-1/p}\, , \]
we obtain for $k(p)$ large enough, thanks to \eqref{eq:equi-integrability},
\[ \dfrac d{dt} Y_p(t) \leq - \delta Y_p(t)^{p/(p-1)} + O\left(Y_p(t)\right) + O(1)\, , \]
where $\delta$ is a positive constant, independant of time.

As a standard consequence, the following estimate holds true for any time $t$ smaller than a reference time $T$:
\[ Y_p(t)\leq C(T) t^{1-p}\, , \]
where the constant $C(T)$ does not depend on the initial value $Y_p(0)$. These {\em a priori} estimates guarantee that the $L^p$ norms of $\rho(t)$ ($p>p_0$) becomes finite for $t>0$.
\end{proof}

\begin{remark}[About the case $1<\alpha\leq 2$]
\label{rem:alpha>1}
The present method can also deal with $1<\alpha\leq 2$, for which global existence has already been proved in \cite{E1}. In fact we shall extend accordingly Proposition \ref{prop:GNS} with $\alpha/2$ derivatives ($\alpha>1$) as following: 
\[
\int_\R \rho^{p+1}(x)\, dx \leq \norm{\rho^{p/2}}{\dot H^{\alpha/2}}^{2\beta} M^{1+p(1-\beta)}\, , \quad \beta = \dfrac p{p+\alpha -1}\, .
\]
Notice  that our strategy requires weaker hypotheses on the initial data 
(in particular regularizing effects can be proved as before).

\end{remark}

\begin{remark}
[Intermediate asymptotics when $\alpha = 1$]
It is known that for the classical two-dimensional Keller-Segel system the cell density in space/time rescaled variables converges to a self-similar profile when mass is subcritical \cite{BDP}. The proof of this fact strongly uses the energy structure. This question is open for the one-dimensional Keller-Segel system with half-diffusion under consideration here.\\ Recall that when only diffusion occurs (without a chemotactic coupling), such a self-similar decay holds true. This can be seen {\em via} the following argumentation in Fourier variables.

First rescale time and space: $u(\tau,y) = (1+t) \rho(t,(1+t)y)$, where $\tau = \log(1+t)$. The new equation reads:
\[ \partial_\tau u(\tau,y) = - \Lambda u(\tau,y) + \partial_y( y u(\tau,y))\, . \]
This writes in Fourier variable as following:
\[ \partial_\tau \hat u(\tau,\xi) = -|\xi| \hat u(\tau,\xi) - \xi \partial_\xi \hat u(\tau,\xi)\, . \]
Or, equivalently,
\[ \partial_\tau\left( \hat u(\tau,\xi) \exp(|\xi|) \right) + \xi \partial_\xi\left(\hat u(\tau,\xi) \exp(|\xi|)\right) = 0\, . \]
As a consequence, $\hat u(\tau,\xi) \exp(|\xi|)$ can be integrated along the characteristics outgoing from 0, where $\hat u(\tau,0)  = M$. This shows that $\hat u(\tau,\xi) \exp(|\xi|)$ converges to $M$ locally in frequency.
Therefore, $u(\tau,y)/M$ converges to the inverse Fourier transform of $\exp(-|\xi|)$, which is nothing but the Cauchy density.

\begin{figure}
\begin{center}
\includegraphics[width = .8\linewidth]{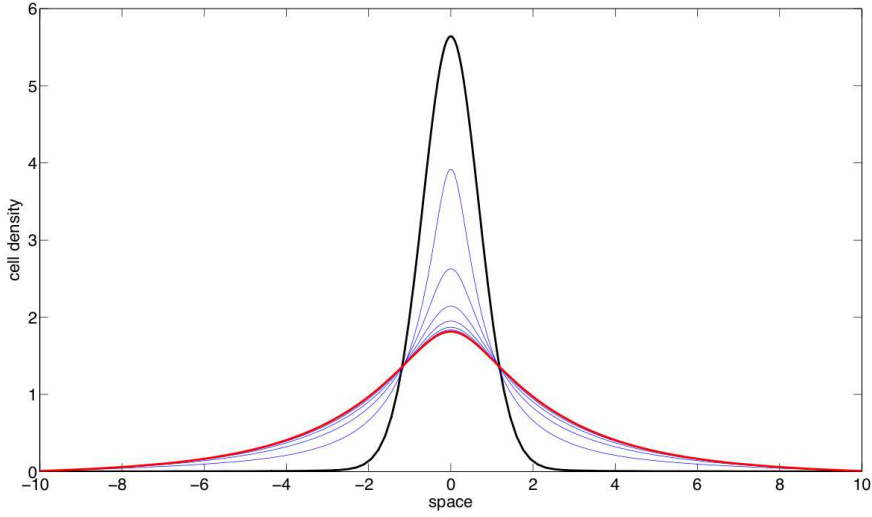}
\caption{Numerical simulation of the one-dimensional Keller Segel system \eqref{frKS} in rescaled variables with $\alpha = 1$. The solution converges to a self-similar profile (in red). Here mass is subcritical as opposed to Fig. \ref{fig:BU}.}
\label{fig:existence}
\end{center}
\end{figure}
Numerical simulations clearly indicate that such a statement is expected to hold true when a chemotactic contribution is added to the diffusion equation and mass is subcritical (see Fig. \ref{fig:existence}).
\end{remark}

\section{Blow-up: proof of Theorem \ref{BU}}

\label{sec:BU}

We focus in this section on the regime $\alpha<1$, for which blow-up may occur. We exhibit a criterion involving 
the mass and the first moment of the initial cell density, in the same spirit as \cite{CPZ}.

Testing the fractional diffusion Keller-Segel against an adequate function $\phi$ (regular with slow decay at infinity to be 
precised below) writes after symmetrization:
\begin{equation}\label{eq:weak formulation}
\begin{split}
\dfrac d{dt}\int_{\R } \phi(x) \rho(t,x)\, dx = \dfrac{c(\alpha)}2 \iint_{\R \times\R } 
\dfrac1{|x-y|^{1+\alpha}}(\phi(x)-\phi(y))(\rho(t,x)-\rho(t,y))\, dxdy\\ - \dfrac14 \iint_{\R \times \R } \sgn(x-y) 
(\phi'(x)-\phi'(y))\rho(t,x)\rho(t,y)\, dxdy\, .   \end{split} \end{equation}

\begin{lemma}[An auxiliary test function] 
\label{lem:auxiliary}
Choose any $0<\beta<1$ satisfying $\alpha+\beta>1$. 
Introduce a $\mathcal C^\infty$, sub-additive, increasing function $\phi$ which satisfies: $\phi(x) = |x|$ for $|x|\leq 1$ and 
$\phi(x) = |x|^{1-\beta}$ for $|x|\geq 2$. Denote $\omega(x) = - \Lambda^\alpha \phi(x)$:
\[ \omega(x) = c(\alpha)\int_{\R } \dfrac{\phi(x+h)-\phi(x)}{|h|^{1+\alpha}}\, dh \, . \]
Then we have the following pointwise estimate for $\omega$: \[\omega(x)\leq C(1+ |x|^{1-\beta})\, .\]
\end{lemma}

\begin{proof}
We split the integral into two parts:
\begin{align*}
|\omega(x)|&\leq \int_{|h|< 2} \dfrac{|\phi(x+h)-\phi(x)|}{|h|^{1+\alpha}}\, dy + \int_{|h|>2} 
\dfrac{|\phi(x+h)-\phi(x)|}{|h|^{1+\alpha}}\, dh  \\
&\leq \int_{|h|< 2} \dfrac{|\phi|_{W^{1,\infty}}}{|h|^{\alpha}}\, dh+ \int_{|h|> 2} \dfrac{\phi(x+h) + \phi(x)}{|h|^{1+\alpha}}\, 
dh \\
& \leq |\phi|_{W^{1,\infty}} C(\alpha) + \int_{|h|> 2} \dfrac{|h|^{1-\beta}+2\phi(x) }{|h|^{1+\alpha}}\, dh \\
& \leq C(|\phi|_{W^{1,\infty}},\alpha) + C(\alpha,\beta) + C(\alpha) \phi(x) \\
& \leq C(|\phi|_{W^{1,\infty}},\alpha,\beta)(1+ |x|^{1-\beta})\, .
\end{align*}

\end{proof}

\begin{proof}[Proof of Theorem \ref{BU}]
The proof begins with testing the Keller-Segel \eqref{eq:weak formulation} against an auxiliary function $\phi_\lambda(x) =  \phi(\lambda x)/\lambda$ where $\phi$ satisfies the assumptions of Lemma \ref{lem:auxiliary}. 
\begin{equation} \label{eq:BU test1}
\dfrac d{dt}\int_{\R} \phi_\lambda(x)\rho(t,x)\, dx  = \int_\R \left( - \Lambda^\alpha\phi_\lambda(x)\right)\rho(t,x)\, dx - \dfrac14 \int_{\R\times \R}  \sgn(x-y)  \left(\phi'_\lambda(x) - \phi'_\lambda(y) \right)\rho(t,x) \rho(t,y)\, dx dy \, .
%
%
\end{equation}
Thanks to a scaling argument and Lemma \ref{lem:auxiliary}, we have the following estimate: 
\[ |\Lambda^\alpha \phi_\lambda(x)| \leq \lambda^{\alpha -1} C\left( 1+ \phi(\lambda x)\right)\, . \] 
As a consequence we have for the first contribution in \eqref{eq:BU test1}:
\[
\abs{\int_\R (- \Lambda^\alpha \phi_\lambda( x)) \rho(t, x) dx} \leq 
C M \lambda^{\alpha - 1} + C \lambda^{\alpha} \int_\R \phi_{\lambda}(x) \rho(t,x)\,  dx\, .
\]


On the other hand, we can write: 
\begin{align*}
& \phi(x) = |x| + R(x)\, , \quad R(x) = \left\{\begin{array}{ll} 0\quad & \mbox{if} \quad |x|<1\\ |x|^{1-\beta} - |x|\, ,\quad & \mbox{if} \quad |x|>2\end{array}\right.\, \\
& \phi'(x) = \sgn(x) + R'(x)\, .
\end{align*}
We have clearly $|R'(x)|\leq C\phi(x)$, hence:
\[
\abs{R'(\lambda x)} \leq C \phi(\lambda x) = C \lambda \phi_{\lambda}(x)\,  . 
\] 
Therefore, we have for the second contribution in \eqref{eq:BU test1}: 
\begin{align*}
- \dfrac14 \int_{\R\times \R}  \sgn(x-y)  \left(\phi'_\lambda(x) - \phi'_\lambda(y) \right)\rho(t,x) \rho(t,y)\, dx dy
& = - \dfrac14 \int_{\R\times \R}  \sgn(x-y)  \left(\sgn(\lambda x) - \sgn(\lambda y) \right)\rho(t,x) \rho(t,y)\, dx dy \\
& \qquad -   \dfrac12 \int_{\R\times \R}  \sgn(x-y)  R'(\lambda x)\rho(x) \rho(y)\, dx dy \\
& \leq - \dfrac12 \int_{\{ (x,y): xy<0 \}}  \rho(t,x) \rho(t,y)\, dx dy +  C M \lambda   \int_{\R}  \phi_\lambda( x)\rho(t,x) \, dx \, . \\ 
\end{align*}
Observe that the symmetry assumption on the cell density $\rho(t,x)$ implies the crucial point:
\begin{align*}   \int_{\{ (x,y): xy<0 \}}  \rho(t,x) \rho(t,y)\, dx dy 
&= 2 \left( \int_{x<0} \rho(t,x)\, dx \right) \left( \int_{y>0} \rho(t,y)\, dy \right) \\
& = \dfrac{M^2}2\, .
\end{align*}
We conclude the above estimates on the `corrected' first moment $I_\lambda(t):= \int \phi_\lambda(x) \rho(t,x) dx $:
\begin{align}
\frac{d I_\lambda}{dt} & \leq  {C M \lambda^{\alpha - 1 }} + C \lambda^{\alpha} I_{\lambda}(t) 
- \frac{M^2}{4} +  C \lambda M I_{\lambda}(t) \nonumber \\
& \leq  \dfrac{M}{ 4\lambda } \left( C   \lambda^{\alpha}  - \lambda M \right)
+ C \left( \lambda^{\alpha} + \lambda M \right) I_{\lambda}(t) \, . \label{P}
\end{align}
We now choose $\lambda$ such that the terms $\lambda^\alpha$ and $\lambda M$ are well-balanced, 
and such that $ C   \lambda^{\alpha}  - \lambda M = - \lambda M/2 $, which  is a negative quantity. 
This leads to $ \lambda = (\mu/M)^{1/(1-\alpha)}$, for some constant $\mu$ depending on $\alpha$, $\beta$, and the specific choice 
of the auxiliary function $\phi$.
Inequality \eqref{P} rewrites:
\[ \frac{d I_\lambda}{dt} \leq -\dfrac{M^2}8 + C \mu^{\alpha/(1-\alpha)} \dfrac{I_\lambda(t)}{M^{\alpha/(1-\alpha)}}\, .  \]
To finish the argumentation, let us observe that imposing a condition of the form 
\begin{equation} \mu^{\alpha/(1-\alpha)} I_\lambda (0)  < C M^{(2-\alpha)/(1-\alpha)} \, , \label{eq:BU crit}\end{equation}
yields that the quantity $I_\lambda$ must vanish in  finite time, which is an obstruction to global existence.

Observe finally that $I_\lambda (0) \leq \int_\R |x|\rho_0(x) dx$, hence \eqref{eq:BU crit} is satisfied if 
$\int_\R |x|\rho_0(x) dx$ is sufficiently small. This completes the proof of Theorem \ref{BU}.

\end{proof}

\begin{figure}
\begin{center}
\includegraphics[width = .8\linewidth]{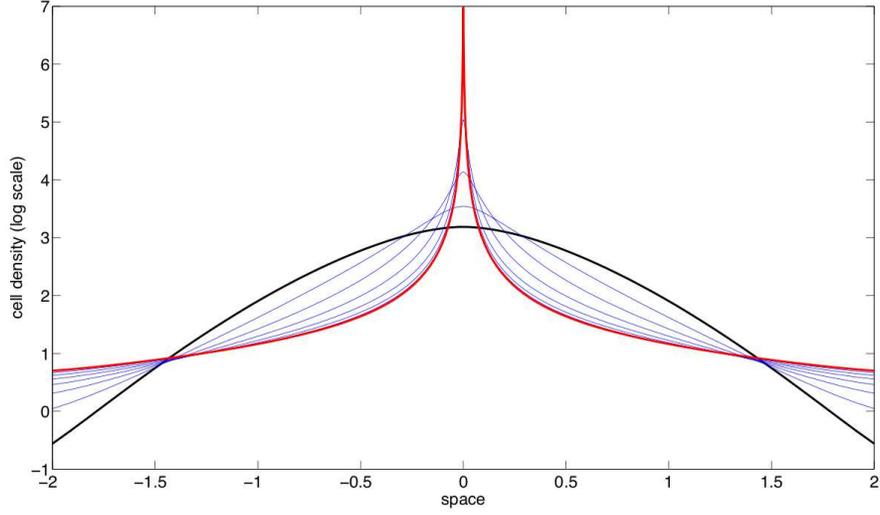}
\caption{Numerical simulation of the one-dimensional Keller Segel system \eqref{frKS} with $\alpha = 1$ for large mass (cell density is plotted in logarithmic scale). The solution clearly blows-up (final profile is plotted in red).}
\label{fig:BU}
\end{center}
\end{figure}

\begin{remark}[On the constants as $\alpha \nearrow 1$]
Tracking carefully the constants in the preceding proof, it turns out that $\mu$ scales like $1/(1-\alpha)$ whereas other constants are indeed of order 1. Thus criterion \eqref{eq:BU crit} rewrites:
\[
I_\lambda (0)^{1-\alpha}  < C^{1-\alpha} (1-\alpha)^{ \alpha }  M^{(2-\alpha) }\, .
\] 
This clearly shows that the previous argument is not expected to be extended to the case $\alpha = 1$. However numerical simulations clearly show that a critical mass is likely to occur when $\alpha = 1$ (see Fig. \ref{fig:BU}).
\end{remark}

\end{document}